\documentclass[12pt]{amsart}
\usepackage{amssymb,latexsym,eufrak,amsmath,amscd, graphicx}

\def\orb{\mathrm{orb}}
\def\st{\mathrm{st}}

\def\<{\langle}
\def\>{\rangle}

\def\ZZ{{\mathbb Z}}
\def\NN{{\mathbb N}}

\def\RR{{\mathbb R}}
\def\QQ{{\mathbb Q}}

\def\cY{\mathcal{Y}}

\def\cT{\mathcal{T}}
\def\cP{\mathcal{P}}

\DeclareMathOperator{\Bbox}{Box}

\DeclareMathOperator{\codim}{codim}

\theoremstyle{plain}
\newtheorem{theorem}{Theorem}[section]

\newtheorem{lemma}[theorem]{Lemma}

\theoremstyle{definition}

\newtheorem{example}[theorem]{Example}

\begin{document}

\title[Ehrhart polynomials and stringy Betti numbers]  {Ehrhart
polynomials and stringy Betti numbers} 

\author[Mustata]{Mircea Musta\c{t}\v{a}}
\address{Department of Mathematics, University of
Michigan\\ Ann Arbor, MI 48109, USA}  
\email{mmustata@umich.edu}

\author[Payne]{Sam Payne}
\address{Department of Mathematics, University of Michigan\\
Ann Arbor, MI 48109, USA}
\email{sdpayne@umich.edu}

\begin{abstract} 
We study the connection between stringy Betti numbers of Gorenstein toric varieties and the generating functions of the Ehrhart polynomials of certain polyhedral regions.  We use this point of view
to give counterexamples to Hibi's conjecture on the
unimodality of $\delta$-vectors of reflexive polytopes.
\end{abstract}

\thanks{The first author was partially supported by 
NSF grant DMS 0500127 and the second author was 
supported by a Graduate Research Fellowship from the NSF}

\maketitle

\section{Introduction}

Let $N$ be a lattice of rank $d$ and let $P$ be a $d$-dimensional
lattice polytope in $N_{\RR} = N \otimes_{\ZZ} \RR$.  For each
nonnegative integer $m$, let $f_{P}(m)$ be the number of lattice
points in $mP$.  Then $f_{P}$ is a polynomial in $m$ of degree $d$,
called the Ehrhart polynomial of $P$.  The generating function
$F_{P}(t) = \sum_{m \geq 0} f_{P}(m) t^{m}$ is a rational function in
$t$ and can be written as
\[
   F_{P}(t) = \frac {\delta_{0} + \delta_{1} t + \cdots + \delta_{d}t^{d}} {
   (1-t)^{d+1} },
\]
for some nonnegative integers $\delta_{i}$, with $\delta_{0} = 1$.  We
put $\delta_{P} = (\delta_{0}, \ldots, \delta_{d})$, and with a slight
abuse of notation we denote by $\delta_P(t)$ the numerator of
$F_P(t)$.  If $\ell$ is the largest $i$ such that $\delta_i$ is
nonzero, then $\ell=d+1-r$, where $r$ is the smallest positive integer
such that $rP$ contains a lattice point in its interior.  Recall that
a lattice polytope is reflexive if it contains $0$ in its interior and
its polar polytope has vertices in the dual lattice.  Given the
lattice polytope $P$, we have $\delta_i=\delta_{\ell-i}$ for all $i$
if and only if $rP$ is the translate of a reflexive polytope.  Hibi
conjectured in \cite[p.\ 111]{Hibi92} that if this is the case, then
$\delta_P$ is unimodal:
\begin{equation} \label{strong inequality}
   \delta_{0} \leq \cdots \leq \delta_{[\ell/2]}.
\end{equation}
In the particular case when $P$ is the Birkhoff polytope of doubly stochastic
$n\times n$ matrices, unimodality had been conjectured by Stanley \cite{St} and
was recently proved by Athanasiadis \cite{At}.

We assume now that $P$ is reflexive, so $\ell=d$.
Hibi showed that in this case
\begin{equation} \label{weak inequality}
   \delta_{0} \leq \delta_{1} \leq \delta_{j}
\end{equation}
for $2 \leq j\leq [d/2]$.  
If, in addition, the boundary of $P$ admits a regular
triangulation such that the vertices of each facet are a basis for the
lattice $N$, then $\delta_{P}$ is the $h$-vector of the triangulation
(see \cite{Hibi}).  In particular, if such a triangulation exists, then
Stanley's theorem on the $h$-vectors of simplicial polytopes
implies that $\delta_{P}$ is unimodal, so $P$ satisfies Hibi's conjecture.

Note that if $P$ is a reflexive polytope of
dimension $d \leq 5$, then Hibi's conjecture follows from (\ref{weak
inequality}).  The following reflexive polytope
gives a counterexample to the conjecture
for $d = 6$; for a more restricted and still open
version of the conjecture, see \cite{Hibi05}.

\begin{example} \label{counterexample}
Let $f = \frac{1}{3}(e_{1} + \cdots + e_{6})$ in $\RR^{6}$, let $N$ be
the lattice $\ZZ^{6} + \ZZ \cdot f$, and let $P$ be the polyope with
vertices $\{ e_{1}, \ldots, e_{6}, e_{1} - f, \ldots, e_{6} - f \}$.
It is straightforward to check that $P$ is reflexive, and one computes
that $2P$ and $3P$ contain 78 lattice points and 314 lattice points,
respectively.  It follows that $\delta_{P} = (1, 6, 8, 6, 8, 6, 1)$.
\end{example}
\noindent In this paper we give a combinatorial formula for
$\delta_{P}$ when $P$ is reflexive, as a positive linear combination of
shifted $h$-vectors of simplicial polytopes, which we arrive at by using toric varieties
to equate the combinatorial invariants $\delta_{i}$ of $P$ with
``stringy'' invariants from complex algebraic geometry.  This formula can also be proved directly, using elementary combinatorial arguments.  We present proofs from both points of view.  With this
formula in hand, it is not difficult to construct examples, such as
Example \ref{counterexample}, where $\delta_{P}$ is not unimodal.

In order to explain our approach, we first reinterpret 
in algebro-geometric terms the proof of unimodality of $\delta_P$
in the special case mentioned above, due to Hibi. 
Here and throughout, $P$
is assumed to be reflexive unless stated otherwise.  
Since $P$ is reflexive, the polar polytope $P^{\circ}$ is reflexive, too.
Note that the polytope $P^{\circ}$ corresponds to a toric variety 
$X_{P^{\circ}}$
defined by the fan over the faces of $P$,
and to an ample divisor $D_{P^{\circ}}$ 
on $X_{P^{\circ}}$. The fact that $P^{\circ}$ is reflexive
is equivalent with the fact that $D_{P^{\circ}}$ is the canonical divisor on
$X_{P^{\circ}}$ 
(so in particular $X_{P^{\circ}}$ is a Fano variety). 

Consider a triangulation ${\mathcal P}$
of  the boundary of
$P$ and let $\Delta$ be the fan whose maximal 
cones are the cones over the facets of ${\mathcal P}$. We have a proper 
birational morphism $f\colon X(\Delta)\to X=X_{P^{\circ}}$
induced by the identity on $N$. If ${\mathcal P}$ is a regular triangulation 
such that the vertices of each facet of ${\mathcal P}$ give a basis of $N$,
then $f$ is a resolution of singularities, $X(\Delta)$ is projective, and $f$
is crepant, 
i.e.\ the pull-back of the canonical bundle on $X$
is isomorphic to the canonical bundle on $X(\Delta)$. Conversely,
every such resolution of singularities of $X$ arises from a triangulation
as above.
Given such a triangulation, $\delta_{i}$ is the $2i$-th
Betti number of $X(\Delta)$, the dimension of the singular
cohomology group $H^{2i}(X(\Delta);\QQ)$, and the unimodality of
$\delta_{P}$ follows from the Hard Lefschetz Theorem on $X(\Delta)$.

In general, there may not exist any crepant resolution of
$X_{P^{\circ}}$.  However, using the theory of motivic integration,
one can define ``stringy Betti numbers'' of $X_{P^{\circ}}$ that agree
with the Betti numbers of a crepant resolution whenever such a
resolution exists \cite{Batyrev}.  A result of Batyrev and Dais shows that $\delta_i$ is the $2i$-th stringy Betti number of $X_{P^\circ}$ \cite[Theorem 7.2]{BD}.  We generalize this result as follows.

If $X=X(\Sigma)$ is a complete, $d$-dimensional Gorenstein toric variety, then there is a function$\Psi_K$ on $N_\RR$ that on each cone is given by an element of the dual lattice, and such that $\Psi_K(v_i) = 1 $ for every primitive generator $v_i$ of a ray of $\Sigma$.  Consider the set
$$Q=\{v\in N_{\RR}\mid \psi_K(v)\leq 1\}.$$  For every cone $\sigma$ in $\Sigma$ the intersection
$Q\cap\sigma$ is a lattice polytope; it is the convex hull of the origin 
and of the primitive generators of the rays of $\sigma$.
We see that $Q$, viewed as the union of the polytopes $Q \cap \sigma$, is naturally a polyhedral  complex, and that $Q$ is homeomorphic to a ball of dimension $d$.

Therefore we may define as in \cite{Hibi95} 
a polynomial of degree $d$ (the Ehrhart polynomial) $f_Q$
such that $f_Q(m)$ is the number of lattice points in $mQ$ for every 
nonnegative integer $m$.
Then we can write the generating function $F_Q(t) = \sum_{m \geq 0} f_Q(m) t^m$ in the form
\[
   F_Q(t) = \frac{ \delta_0  + \delta_1 t + \cdots + \delta_d t^d } { (1-t)^{d+1} },
\]   
for some nonnegative integers $\delta_i$.

\begin{theorem}\label{thm1}
For every complete Gorenstein toric variety $X$, $\delta_i$ is equal to the $2i$-th
stringy Betti number of $X$.
\end{theorem}

Although there may not
exist any crepant resolution of singularities for $X$, we
can always find a projective crepant morphism of toric varieties $f: X(\Delta)
\rightarrow X$ such that $X(\Delta)$ has
only Gorenstein orbifold singularities.  Since $f$ is crepant, the
stringy Betti numbers of $X$ are equal to the stringy
Betti numbers of $X(\Delta)$.  A theorem of Yasuda \cite{Yasuda} then
implies that the stringy Betti numbers of $X(\Delta)$ are equal to the
dimensions of the graded pieces of the orbifold cohomology of
$X(\Delta)$.  We get a combinatorial formula for these dimensions
using a toric formula due to Borisov, Chen, and Smith \cite{BCS}.  The
resulting description of $\delta_Q$ is as follows.  Fix a
triangulation $\cT$ of the boundary of $Q$ whose vertices are in $N$,
and let $\Delta$ be the fan over the faces of $\cT$.  For a face $F
\in \cT$ with vertices $v_{1}, \ldots, v_{r}$, define
\[
   \Bbox (F) = \{ a_{1}v_{1} + \cdots + a_{r}v_{r} \in N_{\RR} : 0 < a_{i}
   < 1 \},
\] 
and let $\Delta_{F}$ be the fan in $N / (N \cap \mbox{span\,} F)$ whose cones
are the projections of the cones in $\Delta$ containing $F$.  For a
positive integer $m$, let $h_{\Delta_{F}}[m]$ denote the $h$-vector of
$\Delta_{F}$ shifted by $m$, defined by
\[
   h_{\Delta_{F}}[m]_{i} = \left\{ \begin{array}{ll} 0 \mbox{ for } i < m. \\
   (h_{\Delta_{F}})_{i-m} \mbox{ for } i \geq m.  \end{array} \right.
\]   
Note that $\Delta_{F}$ is the simplicial fan corresponding to the
$T$-invariant subvariety of $X(\Delta)$ determined by the cone over
$F$, and $(h_{\Delta_{F}})_{i}$ is the $2i$-th Betti number of
$X(\Delta_{F})$.  In particular, if $X(\Delta)$ is projective, then the Hard Lefschetz Theorem on $X(\Delta_{F})$ implies
that $h_{\Delta_{F}}$ is unimodal.

\begin{theorem}\label{thm2}
If $\cT$ is any triangulation of the boundary of $Q$ whose vertices are in $N$, then
\[
   \delta_{Q} = h_{\cT} + \sum_{F \in \cT,\, v \in \Bbox (F) \cap N}
   h_{\Delta_{F}}[\Psi_K(v)].
\]   
\end{theorem}

\noindent In particular, the sum of shifted $h$-vectors in Theorem~\ref{thm2}
is independent of the choice of triangulation.

\section{$\delta$-vectors and stringy Betti numbers}

A $d$-dimensional Gorenstein variety $X$ with canonical singularities 
has a stringy $E$-function
\[
   E_{\st}(X;w,z) \in \ZZ[[w,z]] \cap \QQ(w,z)
\]
defined using Hodge theory and motivic integration on a resolution of
singularities of $X$.  If $E_{\st}(X;w,z) = \sum_{p,q} a_{pq}
w^{p}z^{q}$ is a polynomial, then the $j$-th stringy Betti number of
$X$ is defined to be $(-1)^{j} \sum_{p+q = j} a_{pq}$.  

Suppose now that 
$X=X(\Sigma) $ is
a complete Gorenstein toric variety (see \cite{Fu} for basic facts on toric varieties).
In this case
$E_{\st}(X;w,z)$ is a polynomial in $wz$, so the odd
stringy Betti numbers vanish and the $2i$-th stringy Betti number of
$X$ is the coefficient of $(wz)^{i}$ \cite[Section 3]{Batyrev}.  Our
proof of Theorem \ref{thm1} is based on the following formula for
$E_{\st}(X;w,z)$ as a rational function \cite[Theorem 4.3]{Batyrev}.
Since $X$ is Gorenstein, we have a function $\psi_K$ on $N_\RR$ that on each cone is given by an element in the dual lattice,
and such that
$\Psi_{K}(v_{i}) = 1$ for every primitive generator $v_{i}$ of a ray of
$\Sigma$.  For each cone $\sigma \in \Sigma$, let $\sigma^{\circ}$
denote the relative interior of $\sigma$.  Recall that $\sum_{v \in
\sigma^{\circ}} (wz)^{-\Psi_{K}(v)}$ is a rational function in $wz$
(see, for example, \cite[VIII.1]{Barvinok}).  Batyrev has shown
that we have the following equality of rational functions,
\begin{equation} \label{lattice formula}
    E_{\st}(X;w,z) = (wz - 1)^{d} \sum_{\sigma \in \Delta} \sum_{v \in
    \sigma^{\circ} \cap N} (wz)^{-\Psi_{K}(v)}.
\end{equation}
As in the Introduction, we define
\[
  Q = \{ v \in N_\RR \ | \ \Psi_K(v) \leq 1 \}.
\]
There is an Ehrhart polynomial $f_Q$ such that, for positive integers $m$, $f_Q(m)$ is the number of lattice points in $mQ$, and $f_Q$ satisfies Ehrhart reciprocity: $f_Q(-m)$ is the number of lattice points in the interior of $mQ$.  The proofs of these assertions follow as in \cite{Hibi92}, using the fact that $Q$ is homeomorphic to a ball of dimension $d$.  The generating function $F_Q(t) = \sum _{m \geq 0} f_Q(m) t^m$ can then be written
\[
   F_Q(t)  = \frac { \delta_0 + \delta_1 t  + \cdots + \delta_d t^d } { (1-t)^{d+1} },
\]
for some nonnegative integers $\delta_i$.  

For the proof of Theorem~\ref{thm1} we will need the following two lemmas.  
A proof of the first lemma in the case when $Q$ is a polytope can be found in \cite{Hibi92}
and the general case is similar, but we include the proof
for the reader's convenience.

\begin{lemma}\label{symmetry}
With the above notation, we have $\delta_i=\delta_{d-i}$ for every $i$.
\end{lemma}

\begin{proof}
Note first that if $m$ is a positive integer, then a lattice point $v$ is in
the interior of $mQ$ if and only if $v$ is in $(m-1)Q$. Indeed, $v$ is
in the interior of $mQ$ if and only if $\psi_K(v)<m$, and since 
$\psi_K(v)$ is an integer this is the case if and only if $\psi_K(v)
\leq m-1$, which happens if and only if $v$ is in $(m-1)Q$.
Ehrhart reciprocity implies that
\begin{equation}\label{symmetry1}
f_Q(m-1)=(-1)^df_Q(-m)
\end{equation}
for every positive integer $m$, and therefore for all $m$.

If we write $f_Q(m)=\sum_{i=0}^da_i{{i+m}\choose {i}}$,
then we deduce 
$$F_Q(t)=\sum_{m\in {\mathbb N}}\sum_{i=0}^da_i{{i+m}\choose i}t^m
=\sum_{i=0}^da_i\cdot\sum_{m\in {\mathbb N}}{{i+m}\choose {i}}t^m
=\sum_{i=0}^d\frac{a_i}{(1-t)^{i+1}}.$$
If we put $\widetilde{F}_Q(t)=\sum_{m\geq 1}f_Q(-m)t^m$, then
$$\widetilde{F}_Q(t)=\sum_{i=0}^da_i\cdot \sum_{m\geq i+1}(-1)^i
{{m-1}\choose {i}}t^m
=\sum_{i=0}^d(-1)^i\frac{a_it^{i+1}}{(1-t)^{i+1}},$$ so we have 
the equality of rational functions 
$\widetilde{F}_Q(t)=-F_Q(t^{-1})$. 

On the other hand, (\ref{symmetry1}) gives
$\widetilde{F}_Q(t)=(-1)^dtF_Q(t)$, hence $F_Q(t^{-1})=(-1)^{d+1}tF_Q(t)$.
Since $(1-t)^{d+1}F_Q(t)=\sum_{i=0}^d\delta_it^i$, this equality
gives $\delta_i=\delta_{d-i}$ for every $i$. 
\end{proof}

\begin{lemma} \label{lattice sum}
With the above notation, we have
    \[
       (1-t)F_{Q}(t) = \sum_{v \in N} t^{\psi_K(v)}.
    \]
\end{lemma}

\begin{proof}
    We can write 
$$F_{Q}(t) = \sum_{m \in \NN} \sum_{v \in mQ \cap N}
    t^{m}=\sum_{v\in N}\sum_{m\geq \psi_K(v)}t^m,$$
using the fact that $v$ is in $mQ$ if and only if $m\geq\psi_K(v)$.  
The assertion in the lemma follows.
\end{proof}

\begin{proof}[Proof of Theorem \ref{thm1}.] 
It is enough to show
that $E_{\st}(X;t, 1) = \delta_{Q}(t)$. 
Combining Lemma \ref{lattice sum} with (\ref{lattice formula}), we
have
\[
   E_{\st}(X;t,1) = (t - 1)^{d} (1-t^{-1}) F_{Q}(t^{-1}).
\]
Now
\[
   F_{Q}(t^{-1}) = \frac{\delta_{Q}(t^{-1})}{(1-t^{-1})^{d+1}}.
\]
By Lemma~\ref{symmetry} we have
 $\delta_{i} = \delta_{d-i}$, so $\delta_{Q}(t^{-1}) = t^{-d}
\delta_{Q}(t)$.  Hence
$$
   E_{\st}(X;t,1)=(t-1)^{d} \frac{ \delta_{Q}(t) } {
   t^{d}(1-t^{-1})^{d}}=  \delta_{Q}(t).
$$ 
\end{proof}

\section{$\delta$-vectors via orbifold cohomology}

The orbifold cohomology of a Gorenstein variety $Y$ with quotient
singularities was defined by Chen and Ruan \cite{CR} and
Yasuda \cite{Yasuda}, as follows.
There is a canonically associated
orbifold (smooth Deligne-Mumford stack) $\cY$ whose coarse moduli
space is $Y$.  Let $I(\cY)$ be the inertia stack of $\cY$.  We denote by
$\cY_{i} \subset I(\cY)$  the connected components of $I(\cY)$ and let
$\overline{\cY_{i}}$ be the coarse moduli space of $\cY_{i}$. 
The ``age'' $s_i$ of $\cY_{i}$ is a positive integer
determined by the action of the inertia group.   As a graded
vector space, the orbifold cohomology of $Y$ is given by
\[
   H^{*}_{\orb}(Y,\QQ) = \bigoplus_{\cY_{i} \subset I(\cY)}
   H^{*}(\overline{\cY_{i}},\QQ) [s_{i}],
\]
where $[s_{i}]$ denotes a grading shift by $s_{i}$, so
$H^{j}(\overline{\cY_{i}},
\QQ) [s_{i}] = H^{j-s_{i}}(\overline{\cY_{i}}, \QQ)$.

It is a theorem of Yasuda \cite{Yasuda} that the $j$-th stringy 
Betti number of $Y$ is equal to the dimension of $H^j_{\orb}(Y,\QQ)$.
See also \cite{Poddar} for a proof of this result in the case
of toric varieties. 

We mention that Chen and Ruan have constructed a ring structure 
on orbifold cohomology in \cite{CR}.
J.~Fernandez gave in \cite{Fe} 
a necessary and sufficient condition for when the
Chen-Ruan cohomology satisfies the Hard Lefschetz Theorem. 
His condition inspired us in looking for the counterexamples
to Hibi's Conjecture.

There is an algebraic version of orbifold cohomology, due to 
Abramovich, Graber and Vistoli \cite{AGV}, the so-called 
orbifold Chow ring. Note that when $Y$ is a simplicial
toric variety, each $\overline{\cY_{i}}$ is also a simplicial toric
variety, so the odd cohomology of $\overline{\cY_{i}}$ vanishes and
$H^{2*}(\overline{\cY_{i}}, \QQ)$ is isomorphic to the Chow 
ring $A^{*}(\overline{\cY_{i}}, \QQ)$.  It follows that
at least as vector spaces,
$H^{2*}_{\orb}(Y, \QQ)$ agrees in this case with 
the \cite{AGV} version $A^{*}_{\orb}(Y,\QQ)$
as used by Borisov, Chen and Smith \cite{BCS}.
We mention that 
while there seems
to be agreement among experts that the ring structures are 
also the same in this case,
there is no available reference. We stress however that  we do not need
this, as we are interested only in the vector space structure of
the orbifold cohomology.

\begin{proof}[Proof of Theorem \ref{thm2}.] 
Let $Y$ be the toric variety corresponding to the fan $\Delta$ whose
maximal cones are the cones over the facets of the triangulation
$\cT$.  For a face $F \in \cT$ and a lattice point $v \in \Bbox(F)$,
$\Delta_{F}$ is the fan associated to the stacky fan
$\mathbf{\Delta}/\mathbf{\sigma}(\overline{v})$ defined in
\cite{BCS}, and hence $h_{\Delta_{F}}$ is the vector whose $i$-th
entry is the dimension of $A^{i}(X(\Delta_{F}))$.  Furthermore, the
integer $\Psi_K(v)$ is equal to $\deg y^{v}$ as defined in \cite{BCS}.
Hence the theorem follows from \cite[Proposition 5.2]{BCS}.
\end{proof}

Although we arrived at Theorem \ref{thm2} through the connection with orbifold cohomology and the results of \cite{BCS}, it is also possible to prove this result directly using elementary combinatorial methods, as follows.  For a fan $\Delta$ with $h$-vector $h_\Delta = (h_0, \ldots, h_r)$, we write $h_\Delta(t)$ for the polynomial $h_\Delta(t) = h_0 + h_1 t + \cdots + h_r t^r$.
 
\begin{proof}[Second proof of Theorem \ref{thm2}]
 By Lemma \ref{lattice sum}, it will suffice to show that
\[
   (1-t)^d \cdot \sum_{v \in N} t^{\Psi_K(v)} = \sum_{F \in \cT, v \in \Bbox(F)} t^{\Psi_K
   (v)} \cdot h_{\Delta_F}(t) .
\]
Now each lattice point in the cone over a face $G \in \cT$ can be written uniquely as a nonnegative integer linear combination of the vertices of $G$ plus a fractional part.  Hence any lattice point $v_0$ in the relative interior of this cone can be written uniquely as
\[
   v_0 = v + v_{G | F} + v',
\]
where $v$ is in $\Bbox(F)$ for some face $F \prec G$, $v_{G | F}$ is the sum of the vertices of $G$ that are not in $F$, and $v'$ is a nonnegative integer linear combination vertices of $G$.
Since each lattice point $v \in N$ is in the relative interior of exactly one cone, it follows that
\begin{eqnarray*}
 (1-t)^d \sum_{v \in N} t^{\Psi_K(v)} & = & \sum_{F \in \cT, v \in \Bbox({F} ) } t^{\Psi_K(v)} \cdot \sum_{G \succ F} t^{\dim G - \dim F}  (1-t)^{\codim G}  \\
 & = & \sum_{F \in \cT, v \in \Bbox(F)} t^{\Psi_K(v)} \cdot h_{\Delta_F}(t),
\end{eqnarray*}
  as required.
\end{proof}

\begin{example}
    Let $m$ be a positive integer.  Let $f \in \RR^{2m}$ be the vector
    $f = (\frac{1}{m}, \ldots, \frac{1}{m})$, and let $N$ be the
    lattice $N = \ZZ^{2m} + \ZZ \cdot f$.  We take $P \subset \RR^{2m}$ 
to be
    the polytope with vertices $e_{1}, \ldots, e_{2m}, e_{1} - f,
    \ldots, e_{2m} -f $.  It is straightforward to check that $P$ is
    reflexive.  We will show that 
    \[
       \delta_{P} = (1, 2m, 2m+2, 2m, 2m+2, \ldots, 2m , 2m+2, 2m, 1).
    \]   
    This generalizes Example \ref{counterexample}, and shows that for
    $m > 0$ there are $2m$-dimensional reflexive polytopes with
    $[\frac{m-1}{2}]$ descents in $(\delta_{0}, \delta_{1}, \ldots,
    \delta_{m})$.
    
    We compute $\delta_{P}$ by applying Theorem \ref{thm2} to the
    triangulation $\cP$ of the boundary of $P$ whose facets are $\<
    e_{1}, \ldots, e_{2m} \>, \< e_{1} - f, \ldots, e_{2m} - f \>, \<
    e_{1}, \ldots, \widehat{e_{j}}, \ldots, e_{k}, e_{k} - f, \ldots,
    e_{2m-f} \>$ and $\< e_{1}, \ldots, e_{j}, e_{j} - f, \ldots,
    \widehat{e_{k} - f}, \ldots, e_{2m} - f \>$ for $1 \leq j < k \leq
    2m$.  This triangulation is obtained by ``pulling'' the sequence
    of points $e_{1}, \ldots, e_{2m-1}$.  In particular, $\cP$ is a
    regular triangulation, and hence $h_{\cP}$ is unimodal.  Now $\cP$
    has $4m$ vertices, so $(h_{\cP})_{1} = 2m$, and $\cP$ has $4m^{2} -
    2m + 2$ facets, so $(h_{\cP})_{0} + \cdots + (h_{\cP})_{2m} =
    4m^{2} - 2m + 2$.  It then follows from unimodality and the fact
    that $(h_{\cP})_{0} = (h_{\cP})_{2m} = 1$ that 
    \[
    h_{\cP} = (1, 2m, 2m, \ldots, 2m, 2m, 1).
    \]
    To compute $\delta_{P}$, it remains to compute the contributions
    of the points in $\Bbox(F)$ for the faces $F \in \cP$.  The
    only faces of $\cP$ whose $\Bbox$ is nonempty are $F = \< e_{1},
    \ldots, e_{2m} \>$ and $ F' = \< e_{1} - f, \ldots, e_{2m} - f \>$,
    which contain $\{f, \ldots, (m-1)f \}$ and $\{ -f, \ldots, (1-m)f
    \}$, respectively.  Since $F$ and $F'$ are facets, $\Delta_{F} = 
    \Delta_{F'} = 0$ and $h_{\Delta_{F}} = h_{\Delta_{F'}} =
    1$.  Since $m_{v} = 2k$ for $v = \pm k \cdot f$, it follows that
    \[
         \delta_{P} = h_{\cP} + (0,0,2,0,2, \ldots, 2,0,2,0,0),
    \]
    as required.      
\end{example}

\subsection*{Acknowledgements}
 We are grateful to Bill Fulton for bringing Hibi's conjecture to our
 attention, for his suggestions, 
and for his constant encouragement.  We thank the
 referee for pointing out the connection with \cite{BD}, and for 
suggesting the inclusion of a purely combinatorial proof of Theorem \ref{thm2}.

\providecommand{\bysame}{\leavevmode \hbox \o3em

{\hrulefill}\thinspace}

\end{document}